\documentclass[12pt]{amsart}
\usepackage[mathscr]{eucal}
\usepackage{amssymb,amsmath}
\usepackage{amsfonts}
\usepackage{a4}

\begin{document}
\title{On special rational curves in grassmannians.}

\author[Tomasz Maszczyk]{Tomasz Maszczyk\dag}
\address{Institute of Mathematics\\
Polish Academy of Sciences\\
Sniadeckich 8\newline 00--956 Warszawa, Poland\\
\newline Institute of Mathematics\\
University of Warsaw\\ Banacha 2\newline 02--097 Warszawa, Poland}
\email{t.maszczyk@uw.edu.pl}

\thanks{\dag The author was partially supported by KBN grants 1P03A 036 26 and 115/E-343/SPB/6.PR UE/DIE 50/2005-2008.}

\begin{abstract}
We characterize, among all morphisms
$\mathbb{P}^{1}\rightarrow\mathbb{G}(d, 2d)$, those which are
${\rm GL}_{2d}(\mathbb{C})$-equivalent to the canonical morphism
induced by the Morita equivalence
$\mathbb{C}^{d}\otimes_{\mathbb{C}}-$.
\end{abstract}

\maketitle

\paragraph{\textbf{1. Introduction}} Understanding curves, especially rational, in
Grassmann varieties is very helpful in many problems in Algebraic
Geometry \cite{BerDasWen}, \cite{Kir}, \cite{Osl}, \cite{Pap},
\cite{Str}, Algebraic Topology \cite{Hur}, \cite{HurSan},
\cite{ManMil}, and Control Theory \cite{AgrGam}, \cite{AgrZel},
\cite{Cla}, \cite{Byr}, \cite{Hel1}, \cite{Hel2}, \cite{MarHer},
\cite{Zele}, \cite{Zeli}.

In this article we are interested in the characterization of
curves in Grassmann varieties arising as a result of the following
canonical construction.

The Morita equivalence of $\mathbb{C}$-modules and ${\rm
M}_{d}(\mathbb{C})$-modules (where ${\rm M}_{d}(\mathbb{C})$
denotes the ring of $d\times d$ matrices with complex entries) is
defined by sending a $\mathbb{C}$-module $M$ to
$\mathbb{C}^{d}\otimes_{\mathbb{C}}M$. This identifies
grassmannians of $\mathbb{C}$-sub-modules of the
$\mathbb{C}$-module $\mathbb{C}^{d}$ isomorphic to
$\mathbb{C}^{k}$ with grassmannians of ${\rm
M}_{d}(\mathbb{C})$-sub-modules of the ${\rm
M}_{d}(\mathbb{C})$-module
$\mathbb{C}^{dn}=\mathbb{C}^{d}\otimes_{\mathbb{C}}\mathbb{C}^{n}$
isomorphic to
$\mathbb{C}^{dk}=\mathbb{C}^{d}\otimes_{\mathbb{C}}\mathbb{C}^{k}$
$$\mathbb{G}(k,n)\stackrel{\cong}{\rightarrow}\mathbb{G}_{M_{d}}(dk,dn).$$
The faithfull forgetting functor ${\rm M}_{d}(\mathbb{C})-{\rm
mod}\rightarrow\mathbb{C}-{\rm mod}$ defines a canonical closed
embedding $$\mathbb{G}_{M_{d}}(dk,dn)\hookrightarrow
\mathbb{G}(dk,dn).$$

\vspace{3mm}
\paragraph{\textbf{Definition.}}We call the composition
$\mathbb{G}(k,n)\rightarrow\mathbb{G}(dk,dn)$ of the above
morphisms the \textit{Morita morphism}. For $(k,n)=(1,2)$ the
above composition is a smooth rational curve
$$f_{M}:\mathbb{P}^{1}\rightarrow\mathbb{G}(d, 2d),$$
which we call the \textit{Morita curve}.

\vspace{3mm} Its geometry can be described as follows. Let
\begin{align}
0\rightarrow
S\rightarrow\mathbb{C}^{2d}\otimes_{\mathbb{C}}\mathcal{O}_{\mathbb{G}}\rightarrow
Q\rightarrow 0.\label{seq}
\end{align}
be the tautological short exact sequence on the grassmannian
$\mathbb{G}=\mathbb{G}(d, 2d)$, where $\mathcal{O}_{\mathbb{G}}$
denotes the sheaf of holomorphic functions and $S$ (resp. $Q$)
denotes the locally free sheaf of holomorphic sections of the
tautological sub-bundle (resp. quotient bundle) of the trivial
bundle. The fiber $S_{x}$ of the tautological sub-bundle $S$ at a
point $x\in \mathbb{G}$ corresponding to a subspace of
$\mathbb{C}^{2d}$ can be canonically identified with this
subspace. For $d=1$ the short exact sequence (\ref{seq}) becomes
the twisted Euler sequence
\begin{align}
0\rightarrow
\mathcal{O}_{\mathbb{P}^{1}}(-1)\rightarrow\mathbb{C}^{2}\otimes_{\mathbb{C}}\mathcal{O}_{\mathbb{P}^{1}}\rightarrow
{\rm
det}_{\mathbb{C}}(\mathbb{C}^{2})\otimes_{\mathbb{C}}\mathcal{O}_{\mathbb{P}^{1}}(1)\rightarrow
0,
\end{align}
where ${\rm
det}_{\mathbb{C}}(\mathbb{C}^{2})=\bigwedge^{2}_{\mathbb{C}}\mathbb{C}^{2}$
is a one dimensional space transforming under a linear change of
the base vectors in $\mathbb{C}^{2}$ by the determinant of the
transition matrix, and $\mathcal{O}_{\mathbb{P}^{1}}(1)$ is the
dual of the invertible tautological sheaf
$\mathcal{O}_{\mathbb{P}^{1}}(-1)$. .... The sheaf of holomorphic
1-forms $\Omega^{1}_{\mathbb{G}}$ can be expressed as
$\Omega^{1}_{\mathbb{G}}=S\otimes Q^{\vee}$. Since
$\Omega^{1}_{\mathbb{P}^{1}}={\rm
det}_{\mathbb{C}}(\mathbb{C}^{2})^{-1}\otimes_{\mathbb{C}}\mathcal{O}_{\mathbb{P}^{1}}(-2)$
the differential ${\rm d}f:
f^{*}\Omega^{1}_{\mathbb{G}}\rightarrow
\Omega^{1}_{\mathbb{P}^{1}}$ of any rational curve
$f:\mathbb{P}^{1}\rightarrow\mathbb{G}$ is equivalent to a
morphism of locally free sheaves
\begin{align}
f^{*}S\rightarrow {\rm
det}_{\mathbb{C}}(\mathbb{C}^{2})^{-1}\otimes_{\mathbb{C}}f^{*}Q(-2).\label{diff}
\end{align}
By (\ref{seq}) the determinant of (\ref{diff}) is equivalent to an
element
\begin{align}
\Delta (f)\in {\rm
det}_{\mathbb{C}}(\mathbb{C}^{d})^{2}\otimes_{\mathbb{C}}{\rm
H}^{0}(\mathbb{P}^{1},(f^{*}({\rm det}S)(d))^{-2}).\label{delta}
\end{align}
The Pl\"{u}cker embedding $\mathbb{G}\hookrightarrow
\mathbb{P}^{N}$, $N=(\tiny{\begin{array}{c}
          2d \\
          d
          \end{array}
          })-1$,
is defined by means of the very ample line bundle
$\mathcal{L}=({\rm det}S)^{-1}$ and defines the
\textit{Pl\"{u}cker degree} of a morphism $f$
\begin{align}
{\rm deg}(f):=-{\rm deg}f^{*}{\rm det}S.
\end{align}

One has the Grothendieck splitting
\begin{align}
f^{*}S\stackrel{\cong}{\rightarrow}\bigoplus_{i=1}^{d}\mathcal{O}_{\mathbb{P}^{1}}(-a_{i}).\label{split}
\end{align}
We can assume that
\begin{align}
0\leq a_{1}\leq \ldots \leq a_{d},
\end{align}
because $S^{\vee}$ is globally generated. We will denote by
$\varpi(f)$ the \textit{width of the splitting} (\ref{split}),
i.e.
\begin{align}
\varpi(f)=a_{d}-a_{1}.
\end{align}

 In the case of the Morita curve $f_{M}$ the morphism (\ref{diff}) is an
isomorphism and
$f^{*}S=\mathbb{C}^{d}\otimes_{\mathbb{C}}\mathcal{O}_{\mathbb{P}^{1}}(-1)$.
Therefore
\begin{align}
{\rm deg}(f_{M})  = d,\ \ \Delta (f_{M})  \neq 0,\ \
\varpi(f_{M})= 0.\label{num}
\end{align} Note that every $f$ which is
${\rm GL}_{2d}(\mathbb{C})$-equivalent to $f_{M}$ satisfies
(\ref{num}) as well. We prove the following theorem

\vspace{3mm}
\paragraph{\textbf{Theorem.}} {\em Let $f:\mathbb{P}^{1}\rightarrow\mathbb{G}(d,2d)$ be a rational curve such that
\begin{align}
{\rm deg}(f) = d,\ \ \Delta (f)  \neq 0,\ \ \varpi(f)  \leq 3.
\end{align}
Then $f$ is ${\rm GL}_{2d}(\mathbb{C})$-equivalent to $f_{M}$.}

\vspace{3mm} Note that if ${\rm deg}(f) = d$ then the vector space
in (\ref{delta}) is one dimensional. In the case of $f_{M}$ it can
be canonically identified with $\mathbb{C}$ and then $\Delta
(f_{M})=1$.

In the proof of the theorem we use vanishing of a holomorphic
tensor obtained by means of a matrix valued generalization of the
Schwartz derivative. It is a special case of a noncommutative
generalization of the classical Schwartz derivative, where instead
of complex matrices we can use an arbitrary associative (unital)
algebra (e.g. quaternions or Clifford algebras relating it with
conformal mappings) \cite{Mas}. Independently, the matrix valued
Schwartz derivative appeared in Control Theory, where it arises as
an expression of the curvature of a curve in appropriate
coordinates \cite{AgrGam}, \cite{AgrZel}, \cite{Zele},
\cite{Zeli}. A very general abstract approach to a noncommutative
Schwarz derivative was proposed in \cite{RetSha}.

\vspace{3mm}
\paragraph{\textbf{2. Matrix valued Schwartz derivative.}}

\vspace{3mm}
\paragraph{\textbf{Definition.}} Given a holomorphic  matrix valued function
\begin{align}
\mathbb{C}\supset U\rightarrow M_{d}(\mathbb{C}),\\
x\mapsto y=f(x),
\end{align}
whose derivative $y'$ is an invertible matrix valued function we
define the matrix valued quadratic differential
\begin{align}
\sigma(f):=(((y')^{-1}y'')'-\frac{1}{2}((y')^{-1}y'')^{2}){\rm
d}x^{\otimes 2},
\end{align} which we call the \textit{matrix valued Schwartz derivative} of $f$. \vspace{3mm}

\vspace{3mm}
\paragraph{\textbf{Lemma 1.}} {\em $\sigma$ is invariant under
${\rm GL}_{2}(\mathbb{C})$-transformations}
\begin{align}
x\mapsto & \left(\begin{array}{cc}
               \alpha & \beta \\
               \gamma & \delta
               \end{array}\right)\cdot x=\frac{\alpha x + \beta }{\gamma x+
               \delta},\ \ \ \left(\begin{array}{cc}
               \alpha & \beta \\
               \gamma & \delta
               \end{array}\right)\in {\rm
               GL}_{2}(\mathbb{C}),\label{GL2}
\end{align}
{\em under ${\rm GL}_{2d}(\mathbb{C})$-transformations}
\begin{align}
y\mapsto & \left(\begin{array}{cc}
               A & B \\
               C & D
               \end{array}\right)\cdot y=(Ay+B)(Cy+D)^{-1},\ \ \ \left(\begin{array}{cc}
               A & B \\
               C & D
               \end{array}\right)\in {\rm
               GL}_{2d}(\mathbb{C}),\label{GL2d}
\end{align}
{\em transforms as follows}
\begin{align}
\sigma(\left(\begin{array}{cc}
               A & B \\
               C & D
               \end{array}\right)\cdot
               f)=(Cf+D)\sigma(f)(Cf+D)^{-1},\label{sigma}
\end{align}
 {\em and vanishes if and only if}
\begin{align}
f(x)=(Ax+B)(Cx+D)^{-1},\ \ \ \left(\begin{array}{cc}
               A & B \\
               C & D
               \end{array}\right)\in {\rm
               GL}_{2d}(\mathbb{C}).\label{moebius}
\end{align}

\textit{Proof:} The transformation rules can be checked by
straightforward computation. Vanishing of $\sigma(f)$ is
equivalent to the following system of ODE's

\begin{align}
y'' & = 2y'z,\label{system'}\\
z' & = z^{2}.\label{system''}
\end{align}
On one hand, by (\ref{sigma}) the map of the form (\ref{moebius})
is a solution to the system (\ref{system'})-(\ref{system''}),
because the map $f(x)=x\cdot 1_{d}$ is such. On the other hand,
the solution (\ref{moebius}) with
\begin{align}
\left(\begin{array}{cc}
               A & B \\
               C & D
               \end{array}\right)=\left(\begin{array}{cc}
               y'_{0}- y_{0}z_{0} & y_{0}+y_{0}z_{0}x_{0}-y'_{0}x_{0} \\
               -z_{0} & 1+z_{0}x_{0}
               \end{array}\right)
\end{align}
satisfies the initial condition $(x, y, y', z)=(x_{0}, y_{0},
y'_{0}, z_{0})$. $\Box$

\vspace{3mm}
\paragraph{\textbf{Corollary 1.}} {\em Let $f:\mathbb{P}^{1}\rightarrow\mathbb{G}(d,
2d)$ be a morphism such that}

\begin{align}
{\rm deg}(f) = d,\ \ \Delta (f)  \neq 0.
\end{align}
{\em Then

1) the open subset $U:=f^{-1}(M_{d}(\mathbb{C}))\subset
\mathbb{P}^{1}$, the pre-image of a big affine cell
$M_{d}(\mathbb{C})\subset \mathbb{G}(d, 2d)$, is contained in
$\mathbb{C}\subset \mathbb{P}^{1}$,

2) the first derivative of $f\mid_{U}:U\rightarrow
M_{d}(\mathbb{C})$ is an invertible matrix valued function,

3) $\sigma(f\mid_{U})$ extends uniquely to an element}
\begin{align}
\sigma(f)\in {\rm
det}_{\mathbb{C}}(\mathbb{C}^{2})^{-2}\otimes_{\mathbb{C}}{\rm
H}^{0}(\mathbb{P}^{1},f^{*}(\mathcal{E}nd(S))(-4)).
\end{align}

\textit{Proof:} Since ${\rm deg}(f) = d$ the sheaf ${\rm
det}_{\mathbb{C}}(\mathbb{C}^{d})^{2}\otimes_{\mathbb{C}}(f^{*}({\rm
det}S)(d))^{-2}$ is a trivial line bundle. Therefore, if $\Delta
(f)  \neq 0$ then it is non-zero at every point, hence
(\ref{diff}) is an isomorphism. In particular $f$ is non-constant,
what implies 1).

Every grassmannian $\mathbb{G}(k,n)$ admits an  atlas by big
affine cells $M_{k\times(n-k)}(\mathbb{C})$ with transition
functions of the form
\begin{align}
y\mapsto & \left(\begin{array}{cc}
               A & B \\
               C & D
               \end{array}\right)\cdot y=(Ay+B)(Cy+D)^{-1},\ \ \ \left(\begin{array}{cc}
               A & B \\
               C & D
               \end{array}\right)\in {\rm GL}_{n}(\mathbb{C}),
\end{align}
where transition functions for the tautological bundles $S$ and
$Q$ are of the form
\begin{align}
(y,s) & \mapsto  ((Ay+B)(Cy+D)^{-1}, (Cy+D)s),\\
(y,q) & \mapsto  ((Ay+B)(Cy+D)^{-1}, (A-(Ay+B)(Cy+D)^{-1}C)q).
\end{align}
Therefore, since
\begin{align}
((Ay+B)(Cy+D)^{-1})'=(A-(Ay+B)(Cy+D)^{-1}C)y'(Cy+D)^{-1},
\end{align}
the morphism (2) restricted to $U$ can be identified with the
first derivative of $f\mid_{U}:U\rightarrow M_{d}(\mathbb{C})$, so
the fact that it is an isomorphism implies 2). Finally, the
transformation rules (13)-(15) imply 3). $\Box$

\vspace{3mm}
\paragraph{\textbf{Remark.}}  One can assume above that $y^{T}=y$ (resp.
$y^{T}=-y$) and restrict to transformations of the form
\begin{align}
\left(\begin{array}{cc}
               A & B \\
               C & D
               \end{array}\right)\in {\rm Sp}_{2d}(\mathbb{C})\
               \ \left( {\rm resp.}\   {\rm
               O}_{2d}(\mathbb{C})\right),
\end{align}
to obtain the analogical atlas for lagrangian (resp. isotropic)
grassmannians. One can also use the formula (12) to define the
lagrangian (resp. isotropic) analog of the matrix valued Schwartz
derivative. If we assume that $y$ is a point of the subset
$\mathcal{D}_{d}$ of the lagrangian grassmannian, consisting of
positive definite lagrangian subspaces, i.e. $y^{T}=y$ and
$I_{d}-y\bar{y}$  positive definite, and restrict to the
transformations with
\begin{align}
\left(\begin{array}{cc}
               A & B \\
               C & D
               \end{array}\right)\in \left(\begin{array}{cc}
               iI_{d} & I_{d} \\
               I_{d} & iI_{d}
               \end{array}\right){\rm Sp}_{2d}(\mathbb{Z})\left(\begin{array}{cc}
               iI_{d} & I_{d} \\
               I_{d} & iI_{d}
               \end{array}\right)^{-1}\subset {\rm Sp}_{2d}(\mathbb{C}),
\end{align}
we obtain the global structure of a quotient $\mathcal{A}_{d}={\rm
Sp}_{2d}(\mathbb{Z})\setminus \mathcal{D}_{d}$ parameterizing
isomorphism classes of complex principally polarized abelian
varieties of dimension $d$ \cite{Del}. In all such cases the
restriction of the Euler sequence gives the short exact sequence
\begin{align}
0\rightarrow
S\rightarrow\mathbb{C}^{2d}\otimes_{\mathbb{C}}\mathcal{O}\rightarrow
S^{\vee}\rightarrow 0,
\end{align}
where $S$ is the tautological bundle of lagrangian (resp.
isotropic) subspaces in $\mathbb{C}^{2d}$ with the symplectic
(resp. quadratic form) represented by the standard skew symmetric
(resp. symmetric) matrix
\begin{align}
\left(\begin{array}{cc}
               0 & I_{d} \\
               -I_{d} & 0
               \end{array}\right)\
               \ \left( {\rm resp.}\ \left(\begin{array}{cc}
               0 & I_{d} \\
               I_{d} & 0
               \end{array}\right)\right).
\end{align}

The tangent sheaf is then of the form $\mathcal{T}={\rm
Sym}^{2}(S^{\vee})$ (resp. $\mathcal{T}=\bigwedge^{2}(S^{\vee})$).
By the Uniformization Theorem every Riemann surface $C$ admits an
atlas with transition functions of the form (13). Therefore our
matrix valued Schwarz derivative is well defined for all
holomorphic mappings $f$ from $C$ into usual grassmannians
$\mathbb{G}(d,2d)$, lagrangian or isotropic grassmannians, and
Shimura varieties $\mathcal{A}_{d}$. Then
\begin{align}
\sigma(f)\in {\rm H}^{0}(C,f^{*}(\mathcal{E}nd(S))(2K_{C})).
\end{align}

 \vspace{3mm}
\paragraph{\textbf{Corollary 2.}} {\em Let $f$ be such as in Corollary 1.
If $\sigma(f)=0$ then $f$ is ${\rm
GL}_{2d}(\mathbb{C})$-equivalent to $f_{M}$.}

\vspace{3mm} \textit{Proof:} There exist big affine cells of the
grassmannians $\mathbb{P}^{1}$ and $\mathbb{G}(d,2d)$ for which
the restriction of $f_{M}$ takes the form
\begin{align}
x\mapsto  x\cdot 1_{d}.
\end{align}
Since $\sigma(f)=0$ we know by Lemma 1 that the restriction of $f$
has to be of the form (16), which means that there exists an
element $g\in {\rm GL}_{2d}(\mathbb{C})$ such that $f=g\circ f_{M}
$. $\Box$

\vspace{3mm}
\paragraph{\textbf{3. Proof of the theorem.}} Using the
Grothendieck splitting (5) we see that by (21) $\sigma(f)$ is an
element of a vector space isomorphic to
\begin{align}
\bigoplus_{i,j=0}^{d}{\rm
H}^{0}(\mathbb{P}^{1},\mathcal{O}_{\mathbb{P}^{1}}(a_{i}-a_{j}-4)).
\end{align}
If $\varpi(f)\leq 3$ then $a_{i}-a_{j}-4<0$, for all $i,j$.
Therefore the vector space (21) is zero, hence $\sigma(f)=0$. By
Corollary 2, this implies that $f$ is ${\rm
GL}_{2d}(\mathbb{C})$-equivalent to $f_{M}$. $\Box$

\vspace{3mm}
\paragraph{\textbf{4. Particular cases.}} Since $\varpi(f)\leq d$,
the condition $\varpi(f)\leq 3$ is satisfied automatically for
$d\leq 3$. For $d=1$ or $2$ our theorem reduces to the following
simple and well known facts.

For $d=1$ the Pl\"{u}cker embedding is an identity of
$\mathbb{P}^{1}$, the conditions ${\rm deg}(f)=1$, $\Delta(f)\neq
0$ mean that $f$ is birational \'{e}tale. The Morita curve $f_{M}$
is the identity. Then our theorem is equivalent to the fact that
$f$ is a M\"{o}bius map.

For $d=2$ the  Pl\"{u}cker embedding $\mathbb{G}(2,4)\hookrightarrow
\mathbb{P}^{5}
=\mathbb{P}(\bigwedge^{2}(\mathbb{C}^{2}\otimes_{\mathbb{C}}\mathbb{C}^{2})$
identifies the grassmannian with the Klein quadric
\begin{align}
z_{11,12}z_{21,22}-z_{11,21}z_{12,22}+z_{11,22}z_{12,21}=0.
\end{align}
Then the image of the Morita curve is a closed embedding onto the
intersection of the Klein quadric with the plane
\begin{align}
z_{11,21}=0,\ \ z_{12,22}=0,\ \ z_{11,22}+z_{12,21}=0.
\end{align}
The conditions ${\rm deg}(f)=2$, $\Delta(f)\neq 0$ are equivalent
to the condition that $f$ is a closed embedding onto a conic not
contained in a plane entirely contained in the Klein quadric. Then
our theorem is equivalent to the fact that all conics in
$\mathbb{P}^{5}$ lying on the smooth quadric but not contained in
a plane entirely contained in the quadric, are equivalent up to
automorphisms of $\mathbb{P}^{5}$ preserving the quadric.

For $d=3$ our theorem says about cubics lying on some
9-dimensional smooth intersection of quadrics in
$\mathbb{P}^{19}$, but seems not to reduce to any simple classical
fact.

\end{document}